\documentclass[11pt]{amsart}
\usepackage{amsmath,amssymb}
\usepackage{graphicx}

\usepackage{fullpage}
\usepackage{enumerate}
\usepackage{pgf,tikz}
\usetikzlibrary{arrows}
\usetikzlibrary{positioning}
\usetikzlibrary{patterns}

\def\COMMENT#1{}
\let\COMMENT=\footnote
\def\TASK#1{}

\newdimen\margin   
\def\textno#1&#2\par{%
    \margin=\hsize
    \advance\margin by -4\parindent
           \setbox1=\hbox{\sl#1}%
    \ifdim\wd1 < \margin
       $$\box1\eqno#2$$%
    \else
       \bigbreak
       \hbox to \hsize{\indent$\vcenter{\advance\hsize by -3\parindent
       \sl\noindent#1}\hfil#2$}%
       \bigbreak
    \fi}

\newtheorem{thm}{Theorem}[section]
\newtheorem{define}[thm]{Definition}

\newtheorem{lem}[thm]{Lemma}

\newtheorem{fact}[thm]{Fact}

\newtheorem{question}[thm]{Question}

\newtheorem*{thm*}{Theorem}
\newtheorem*{define*}{Definition}
\newtheorem*{examp*}{Example}
\newtheorem*{lem*}{Lemma}
\newtheorem*{claim*}{Claim}
\newtheorem*{fact*}{Fact}
\newtheorem*{col*}{Corollary}
\newtheorem*{conj*}{Conjecture}
\newtheorem{exex}{Extremal Example}

\begin{document}

\title{A note on colour-bias Hamilton cycles in dense graphs}

\author{Andrea Freschi, Joseph Hyde, Joanna Lada and Andrew Treglown}

\date{}
\begin{abstract}
Balogh, Csaba, Jing and Pluh\'ar recently determined the minimum degree threshold that ensures a $2$-coloured graph $G$ contains a Hamilton cycle of significant colour bias (i.e., a Hamilton cycle that contains significantly more than half of its edges in one colour). In this short note we extend this result, determining the corresponding threshold for $r$-colourings.
\end{abstract}
\maketitle

\section{Introduction} \label{Introduction}
The study of colour-biased structures in graphs concerns the following problem.
Given graphs $H$ and $G$, what is the largest $t$ such that in any $r$-colouring of the edges of $G$, there is always a copy of $H$ in $G$ that has at least $t$ edges of the same colour? Note if $H$ is a subgraph of $G$, one can trivially ensure a copy of $H$ with at least $|E(H)|/r$ edges of the same colour; so one is  interested in when one can achieve a colour-bias significantly above this.

The topic was first raised by Erd\H{o}s in the 1960s (see~\cite{e2, e1}).
  Erd\H{o}s, F\"uredi, Loebl and S\'os~\cite{EFLS}  proved the following:
 for some constant $c>0$, given any $2$-colouring  of the edges of $K_n$ and any 
 fixed spanning tree $T_n$ with maximum degree $\Delta$, $K_n$ contains a copy of $T_n$ such that at least $(n-1)/2+c(n-1-\Delta)$ edges of this copy of $T_n$ receive the same colour. In~\cite{BCsJP}, Balogh, Csaba, Jing and Pluh\'ar
investigated the colour-bias problem  in the case of 
 spanning trees, paths and Hamilton cycles
for various classes of graphs $G.$ Note all their results concern $2$-colourings and therefore were expressed in the equivalent language of \emph{graph discrepancy}.
The following result determines the minimum degree threshold
for forcing a Hamilton cycle of significant colour bias in a 2-edge-coloured graph.

\begin{thm}[Balogh, Csaba, Jing and Pluh\'ar~\cite{BCsJP}]\label{old}
Let  $0<c<1/4$  and $n\in \mathbb N$ be sufficiently large. If $G$ is an $n$-vertex graph with $$\delta(G)\geq(3/4+c)n,$$
 then given any $2$-colouring of $E(G)$ 
there is a Hamilton cycle in $G$ with at least 
$ (1/2+c/64)n$ edges of the same colour.
Moreover, if $4$ divides $n$, there is an $n$-vertex graph $G'$ with $\delta(G')=3n/4$ and
a $2$-colouring of $E(G')$ for which every Hamilton cycle in $G'$ has precisely $n/2$ edges in each colour.
\end{thm}
In~\cite{newkri}, Gishboliner,  Krivelevich and  Michaeli considered colour-bias Hamilton cycles in the random graph $G(n,p)$. Roughly speaking, their result states that if $p$ is such that with high probability (w.h.p.) $G(n,p)$ has a Hamilton cycle, then in fact w.h.p., given any $r$-colouring of the edges of $G(n,p)$, one can guarantee a Hamilton cycle that is essentially as colour-bias as possible (see~\cite[Theorem 1.1]{newkri} for the precise statement). A discrepancy (therefore colour-bias) version of the Hajnal--Szemer\'edi theorem was proven in~\cite{cpc}.

In this paper we give a very short proof of the following multicolour generalisation of Theorem~\ref{old}. We require the following definition to state it.

\begin{define}
Let $t,r \in \mathbb{N}$ and $H$ be a graph. We say that an $r$-colouring of the edges of $H$ is \emph{$t$-unbalanced} if at least $|E(H)|/r + t$ edges are coloured with the same colour.
\end{define}

\begin{thm}\label{mainthm}
Let $n, r, d \in \mathbb{N}$ with $r \geq 2$. Let $G$ be an $n$-vertex graph with $\delta(G)\geq \left(\frac{1}{2}+\frac{1}{2r}\right)n + 6dr^2$. Then for every $r$-colouring of  $E(G)$ there exists a $d$-unbalanced Hamilton cycle in $G$. 
\end{thm}
Note  that $n$, $r$ and $d$ may all be comparable in size.  Further, Theorem~\ref{mainthm} implies Theorem~\ref{old} with a slightly better bound on the colour-bias. In the following section we give constructions that show Theorem~\ref{mainthm} is best possible; that is, there are $n$-vertex graphs $G$ with minimum degree $\delta(G)=(1/2+1/2r)n$ such that for some $r$-colouring of $E(G)$, every Hamilton cycle in $G$ uses precisely $n/r$ edges of each colour. The proof of Theorem~\ref{mainthm} is constructive, producing the $d$-unbalanced Hamilton cycle in time polynomial in $n$.

{\bf Remark:} After making our manuscript available online, we learnt of  simultaneous and independent work of Gishboliner,  Krivelevich and  Michaeli~\cite{newkri2}. They prove an asymptotic version of Theorem~\ref{mainthm} (i.e., for sufficiently large graphs $G$) via Szemer\'edi's regularity lemma. They also generalise a number of the results from~\cite{BCsJP}.

\section{The extremal constructions}
Our first extremal example is a generalisation of a $2$-colour construction from~\cite{BCsJP}.
\begin{exex}\label{exex:general}
Let $r, n \in \mathbb{N}$ where $r\geq 2$ and such that $2r$ divides $n$. Then there exists a graph $G$ on $n$ vertices with $\delta(G) = (\frac{1}{2} + \frac{1}{2r})n$, and an $r$-colouring of $E(G)$, such that every Hamilton cycle uses precisely $n/r$ edges of each colour.    
\end{exex}

\begin{proof}
The vertex set of $G$ is partitioned into $r$ sets $V_1, \ldots, V_{r}$ such that $|V_1|=\ldots=|V_{r-1}|=n/2r$, and $|V_r| = (r+1)n/2r$; the edge set of $G$ consists of all edges with at least one endpoint in $V_r$.
Now colour the edges of $G$ with colours $1, \ldots, r$ as follows:
\begin{itemize}
    \item For each $i \in [r-1]$, colour every edge with one endpoint in $V_i$ and one endpoint in $V_{r}$ with colour $i$.
    \item Colour every edge with both endpoints in $V_r$ with colour $r$ 
		(see Figure~1).
\end{itemize}

Observe that $\delta(G) = \left(\frac{1}{2} + \frac{1}{2r}\right)n$, which is 
attained by every vertex in $V_1 \cup \ldots \cup V_{r-1}$. For each $i \in [r-1]$, every vertex in $V_i$ is only adjacent to edges of colour $i$, $|V_i| 
= n/2r$ and $E(G[V_1\cup\ldots\cup V_{r-1}]) = \emptyset$. Hence every Hamilton cycle in $G$ must contain precisely $n/r$ edges of each colour $i \in [r-1]$. Since a Hamilton cycle has $n$ edges, every Hamilton cycle in $G$ must also contain $n/r$ edges of colour $r$. Thus every Hamilton cycle in $G$ uses precisely $n/r$ edges of each colour. 
\end{proof}

\definecolor{qqqqff}{rgb}{0,0,0}
\definecolor{qqzzqq}{rgb}{0.1,0.1,0.1}
\definecolor{ffqqqq}{rgb}{0.8,0.8,0.8}
\definecolor{ududff}{rgb}{0.4,0.4,0.4}
\begin{figure}[!ht]\label{fig1}
\begin{center}
\begin{tikzpicture}[line cap=round,line join=round,>=triangle 45,x=1cm,y=1cm]
\fill[line width=2pt,color=ffqqqq,fill=ffqqqq,fill opacity=1] (-2,6) -- (-2,-2) -- (0,-2) -- (0,6) -- cycle;
\node at (-1,6.5) {$V_3$};
\node at (-1,2) {\LARGE{$\frac{2n}{3}$}};
\fill[line width=2pt, fill opacity=0] (-6,3) -- (-6,1) -- (-4,1) -- (-4,3) -- cycle;
\node at (-5,3.5) {$V_1$};
\node at (-5,2) {\LARGE{$\frac{n}{6}$}};
\fill[line width=2pt, fill opacity=0] (2,3) -- (2,1) -- (4,1) -- (4,3) -- cycle;
\node at (3,3.5) {$V_2$};
\node at (3,2) {\LARGE{$\frac{n}{6}$}};
\fill[line width=1.5pt,color=qqzzqq,fill=qqzzqq,pattern=north west lines,pattern color=qqzzqq] (-4,3) -- (-2,6) -- (-2,-2) -- (-4,1) -- cycle;
\fill[line width=0pt,color=qqqqff,fill=qqqqff,pattern=grid,pattern color=qqqqff] (0,6) -- (2,3) -- (2,1) -- (0,-2) -- cycle;
\draw [line width=2pt] (-2,6)-- (-2,-2);
\draw [line width=2pt] (-2,-2)-- (0,-2);
\draw [line width=2pt] (0,-2)-- (0,6);
\draw [line width=2pt] (0,6)-- (-2,6);
\draw [line width=2pt] (-6,3)-- (-6,1);
\draw [line width=2pt] (-6,1)-- (-4,1);
\draw [line width=2pt] (-4,1)-- (-4,3);
\draw [line width=2pt] (-4,3)-- (-6,3);
\draw [line width=2pt] (2,3)-- (2,1);
\draw [line width=2pt] (2,1)-- (4,1);
\draw [line width=2pt] (4,1)-- (4,3);
\draw [line width=2pt] (4,3)-- (2,3);

\end{tikzpicture}
\caption{Extremal Example \ref{exex:general} for $r = 3$}
\end{center}
\end{figure}
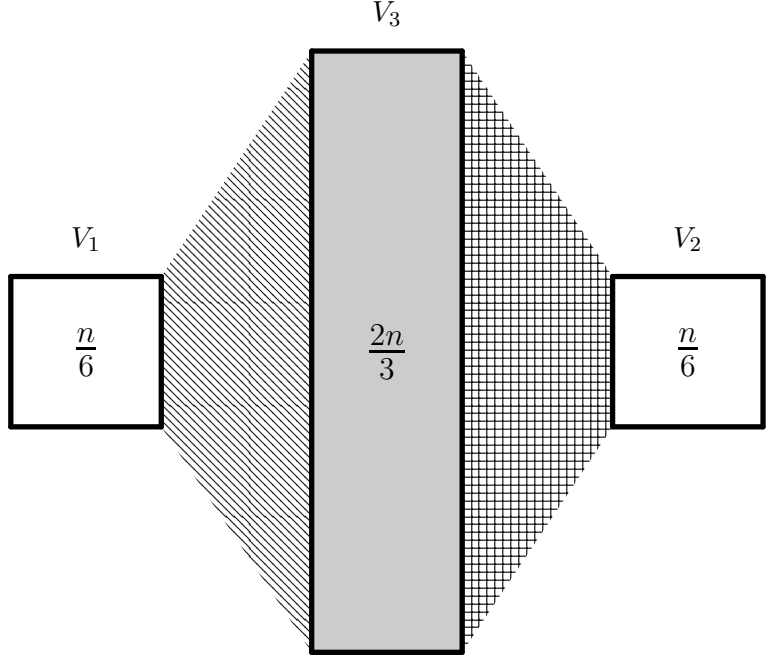

We also have an additional extremal example in the $r=3$ case.

\begin{exex}\label{exex:r=3}
Let $n \in \mathbb{N}$ such that $3$ divides $n$. Then there exists a graph $G$ on $n$ vertices with $\delta(G) = 2n/3$, and a $3$-colouring of $E(G)$, such that every Hamilton cycle uses precisely $n/3$ edges of each colour and every vertex in $G$ is incident to precisely two colours.    
\end{exex}

\begin{proof}
Let $G$ be the $n$-vertex $3$-partite Tur\'an graph. So  $G$  consists of three vertex sets $V_1$, $V_2$ and $V_3$, such that $|V_1| = |V_2| = |V_3| = n/3$, and all possible edges that go between distinct $V_i$ and $V_j$.
Colour all edges between $V_1$ and $V_2$ red; all edges between $V_2$ and $V_3$ blue; all edges between $V_3$ and $V_1$ green.

Clearly $\delta(G) = 2n/3$ and every vertex is incident to precisely two colours. Let $H$ be a Hamilton cycle in $G$ and let $r$, $b$ and $g$ be the number of red, blue and green edges in $H$, respectively. 
Since all red  and green edges in $H$ are incident to vertices in $V_1$, $|V_1| = n/3$ and $V_1$ is an independent set, we must have that $2n/3 = r + g$. Applying similar reasoning to $V_2$ and $V_3$, we have that $2n/3 = b + r$ and $2n/3 = g + b$. Hence $r = b = g = n/3$. Thus every Hamilton cycle in $G$ uses precisely $n/3$ edges of each colour. 
\end{proof}

\section{Proof of Theorem~\ref{mainthm}}
As in~\cite{BCsJP}, 
we  require the following generalisation of Dirac's theorem.

\begin{lem}[P\'{o}sa \cite{p}]\label{lem:forestham}
Let $1 \leq t \leq n/2$, $G$ be an $n$-vertex graph with $\delta(G)\geq \frac{n}{2} + t$ and $E'$ be a set of edges of a linear forest in $G$ with $|E'| \leq  2t$. Then there is a Hamilton cycle in $G$ containing $E'$.
\end{lem}

\begin{proofofmainthm}
Recall that $G$ is a graph on $n$ vertices with $\delta(G) \geq \left(\frac{1}{2}+\frac{1}{2r}\right)n + 6dr^2$ for some integers $r \geq 2$ and $d \geq 1$.
Consider any $r$-colouring of $E(G)$. 
Given a colour $c$ we define the function $L_c:E(G)\rightarrow \{0,1\}$ as follows:
$$L_c(e):=
\begin{cases}
1 \quad \text{if $e$ is coloured with $c$},
\\0 \quad \text{otherwise}.
\end{cases}$$
Given a triangle $xyz$ and a colour $c$, we define $\mbox{Net}_c(xyz,xy)$ as follows:
$$\mbox{Net}_c(xyz,xy) := L_c(xz) + L_c(yz) - L_c(xy).$$ This quantity comes from an operation we will perform later where we extend a cycle $H$ by a vertex $z$ via deleting the edge $xy$ from $H$ and adding the edges $xz$ and $yz$, to form a new cycle $H'$.
One can see that $\mbox{Net}_c(xyz,xy)$ is the change in the number of edges of colour $c$ from $H$ to $H'$. 

Since $\delta(G) \geq \frac{1}{2}n$, by Dirac's theorem, $G$ contains a Hamilton cycle $C$. If $C$ is $d$-unbalanced we are done, so suppose it is not. 
Let $v \in V(G)$. Since $d(v) \geq \left(\frac{1}{2}+\frac{1}{2r}\right)n + 6dr^2$, there are at least $\frac{n}{r}+12dr^2$ edges
$e$ in $C$ such that $v$ and $e$ span a triangle. 

This can be seen in the following way. Let $X$ be the set of neighbours of $v$, and $X^+$ the set of vertices whose `predecessors' on $C$ are neighbours of $v$, having arbitrarily chosen an orientation for $C$. We have 
$$n \geq |X\cup X^+| = |X| + |X^+| - |X\cap X^+| 
    \geq  n + \frac{n}{r} + 12dr^2 - |X\cap X^+|.$$
Hence $|X\cap X^+| \geq \frac{n}{r} + 12dr^2$. Clearly  each element in $X\cap X^+$ yields  a triangle containing $v$, thus giving the desired bound.

This property, together with the fact that $C$ is not $d$-unbalanced (so contains fewer than $n/r+d$ edges of each colour) immediately implies the following.
\begin{fact}\label{facty}
Let $v \in V(G)$, $Y\subseteq V(G)$ with $|Y|\leq 5dr^2$, 
and $xy$ be any edge in $G$ that forms a triangle with $v$ and is disjoint to $Y$.\footnote{Note sometimes in an application of this fact, $xy$ will be an edge of $C$, but other times not.}
Then there is an edge $zw$ on $C$ vertex-disjoint to $xy$, and distinct colours 
$c_{1}$ and $c_{2}$ such that  $vzw$ induces a triangle; $xy$ has colour $c_{1}$; $zw$ has colour $c_{2}$; $z,w \not \in Y$.
\end{fact}

Initially set
 $A := \emptyset$. Consider an arbitrary $v \in V(G)$ and let $x,y,z,w, c_1,c_2$ be as in Fact~\ref{facty} (where $Y:=\emptyset$), where $xy$ is chosen to be an edge of $C$ that forms a triangle with $v$.

If there exists a colour $c$ such that $\mbox{Net}_c(vxy,xy) \not= \mbox {Net}_c(vzw,zw)$ then add the pair $(xy,zw)$ to the set $A$, and define $v_1:=v$. If there is no such colour then we must have that $\mbox{Net}_{c_1}(vxy,xy)=\mbox{Net}_{c_1}(vzw,zw)$ and so
$$L_{c_1}(vx) + L_{c_1}(vy) - L_{c_1}(xy) = L_{c_1}(vw) + L_{c_1}(vz) - L_{c_1}(wz),$$
$$L_{c_1}(vx) + L_{c_1}(vy) - 1 = L_{c_1}(vw) + L_{c_1}(vz)\geq 0,$$ 
as $xy$ has colour $c_1$, $wz$ has colour $c_2$ and $c_1 \neq c_2$.
Hence  $vx$ or $vy$ is coloured with $c_1$. Without loss of generality, let $vx$ be coloured with $c_1$.
By the same argument with colour $c_2$, we may assume that, without loss of generality, $vw$ is coloured $c_2$. Let $c_3$ be the colour of $vy$. Then $\mbox{Net}_{c_3}(vxy,xy) = \mbox{Net}_{c_3}(vzw,zw)$ and so
$$L_{c_3}(vx) + L_{c_3}(vy) - L_{c_3}(xy) = L_{c_3}(vw) + L_{c_3}(vz) - L_{c_3}(wz),$$
$$1 = L_{c_3}(vz),$$ 
as $vx$ and $xy$ are both coloured with $c_1$ and $vw$ and $wz$ are both coloured with $c_2$.
Hence $c_3$ is also the colour of $vz$ (see Figure~2). Since $c_1 \neq c_2$, we may assume, without loss of generality, $c_1\not=c_3$. 

Now we apply Fact~\ref{facty} with $x$ playing the role of $v$; $vy$ playing the role of $xy$; $Y=\emptyset$. We thus obtain a colour $c_4 \not =c_3$ and an edge
$w'z'$ on $C$ that is vertex-disjoint from $vy$, so that $w'z'$ forms a triangle with $x$, and $w'z'$ is coloured $c_4$. Note that by construction  $\mbox{Net}_{c_3}(xvy,vy)=-1$ whilst, as 
$c_4 \not =c_3$, by definition 
$\mbox{Net}_{c_3}(xw'z',w'z')= L_{c_3} (xw')+L_{c_3}(xz')-0\geq 0$.
In this case we define $v_1:=x$ and add 
 the pair $(vy,w'z')$ to $A$.

\definecolor{qqqqff}{rgb}{0.25,0.25,0.25}
\definecolor{ffcctt}{rgb}{0.5,0.5,0.5}
\definecolor{ffqqqq}{rgb}{0.8,0.8,0.8}
\definecolor{yqyqyq}{rgb}{0.5019607843137255,0.5019607843137255,0.5019607843137255}
\definecolor{uququq}{rgb}{0.25098039215686274,0.25098039215686274,0.25098039215686274}
\definecolor{ududff}{rgb}{0.1019607843137255,0.1019607843137255,0.1019607843137255}
\begin{figure}[!ht]\label{c3figure}
\begin{center}
\begin{tikzpicture}[scale = .5, line cap=round,line join=round,>=triangle 45,x=1cm,y=1cm]
\draw [line width=2pt,color=ududff] (0,0) circle (6cm);
\draw [line width=2pt,color=ffqqqq] (0,6)-- (-4.336291901416994,-4.146875033770055);
\draw [line width=2pt,color=ffcctt] (0,6)-- (-2.2142787553979884,-5.576465690147581);
\draw [line width=2pt,color=ffcctt] (0,6)-- (2.190703198420197,-5.585769373723867);
\draw [line width=2pt,color=qqqqff] (0,6)-- (4.341101198561697,-4.141840217083006);
\draw (-0.42377464280574756,7.205627132064676) node[anchor=north west] {$v$};
\draw (-5.25184211764568,-4.158862308645072) node[anchor=north west] {$x$};
\draw (-2.789372130772672,-5.7272167996749605) node[anchor=north west] {$y$};
\draw (1.811946344504827,-5.800591299456177) node[anchor=north west] {$w$};
\draw (4.427667331815402,-4.092236808426289) node[anchor=north west] {$z$};
\draw [shift={(0,0)},line width=2pt,color=ffqqqq]  plot[domain=3.904665970656774:4.334406131420695,variable=\t]({1*6*cos(\t r)+0*6*sin(\t r)},{0*6*cos(\t r)+1*6*sin(\t r)});
\draw [shift={(0,0)},line width=2pt,color=qqqqff]  plot[domain=5.08614767141132:5.52127243451429,variable=\t]({1*6*cos(\t r)+0*6*sin(\t r)},{0*6*cos(\t r)+1*6*sin(\t r)});
\draw [color=ffcctt](-3.9771091250843085,0.5789581675075569) node[anchor=north west] {$c_1$};
\draw [color=ffcctt](-4.164105123334043,-5.203101303831842) node[anchor=north west] {$c_1$};
\draw [color=ffcctt](-1.4580136379925172,-1.9367628198030165) node[anchor=north west] {$c_3$};
\draw [color=ffcctt](0.38,-1.9367373195842333) node[anchor=north west] {$c_3$};
\draw [color=qqqqff](2.709312840785513,0.5789581675075569) node[anchor=north west] {$c_2$};
\draw [color=qqqqff](3.0065558377225478,-5.1464758036130585) node[anchor=north west] {$c_2$};
\draw [ color=uququq](-6.159332113270018,4.9)node[anchor=north west] {\Large${C}$};
\begin{scriptsize}
\draw [fill=uququq] (0,6) circle (2.5pt);
\draw [fill=uququq] (-4.336291901416994,-4.146875033770055) circle (2.5pt);
\draw [fill=uququq] (-2.2142787553979884,-5.576465690147581) circle (2.5pt);
\draw [fill=uququq] (2.190703198420197,-5.585769373723867) circle (2.5pt);
\draw [fill=uququq] (4.341101198561697,-4.141840217083006) circle (2.5pt);
\end{scriptsize}
\end{tikzpicture}
\caption{A Hamilton cycle $C$ for $G$, with a vertex $v$ which is good for $C$. There is no colour $c$ with $\mbox{Net}_c(vxy,xy) \not= \mbox {Net}_c(vzw,zw)$ implying the colour arrangement above. }
\end{center}
\end{figure}
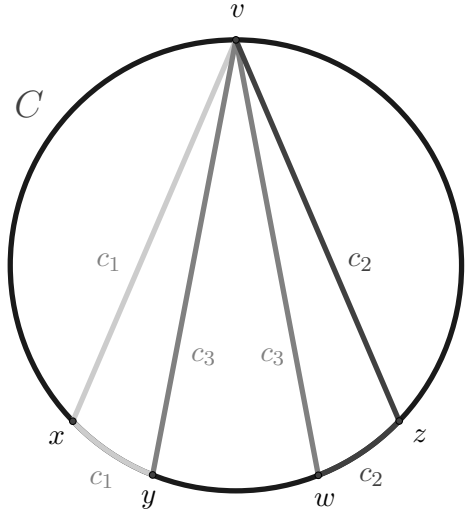

Repeated applications of this argument thus yield
 sets $B:=\{v_1,v_2, \ldots ,v_{dr^2}\}$ and a set  $A$ whose elements are pairs of edges from $G$ so that:
\begin{itemize}
\item All vertices lying in  $B$ and in edges in pairs from $A$ 
are vertex-disjoint.
\item For each $u=v_i$ in $B$ there is a pair $(xy,zw) \in A$ associated with $u$, and a colour $c_u$ so that (i) $uxy$ and $uzw$ are triangles in $G$;  (ii) $\mbox{Net}_{c_u}(uxy,xy) \not= \mbox {Net}_{c_u}(uzw,zw)$. We call $c_u$ the colour \emph{associated} with $u$.
\end{itemize}
Note that it is for the first of these two conditions that we require the set $Y$ in Fact~\ref{facty}. At a given step of our argument, $Y$ will be the set of vertices that have previously been added to $B$ or lie in an edge previously selected for inclusion in a pair from $A$.

\smallskip

There is some colour $c^*$ for which $c^*$ is the colour associated with (at least) $dr$ of the vertices in $B$. Let $B'$ denote the set of such vertices of $B$; without loss of generality we may assume $B'=\{v_1,v_2, \ldots ,v_{dr}\}$. Let $A'$ denote the subset of $A$ that corresponds to $B'$. For each $i \in [dr]$, let $(x_iy_i,z_iw_i)$ 
denote the element of $A'$ associated with $v_i$. We may assume that for each 
 $i \in [dr]$, 
\begin{align}\label{vital}
\mbox{Net}_{c^*}(v_ix_iy_i,x_iy_i) > \mbox {Net}_{c^*}(v_iz_iw_i,z_iw_i).
\end{align}

Consider the induced subgraph $G'$ of $G$ obtained from $G$ by removing the vertices from $B'$. Let $E'$ be the set of all edges which appear in some pair in $A'$. 
As $\delta(G') \geq n/2+dr$,  Lemma~\ref{lem:forestham} implies that there exists a Hamilton cycle $C'$ in $G'$ which contains $E'$.
Let $C_1$ be the Hamilton cycle of $G$ obtained from $C'$ by inserting each $v_i$ from $B'$ between $x_i$ and $y_i$; let $C_2$ be the Hamilton cycle of $G$ obtained from $C'$ by inserting each $v_i$ from $B'$ between $z_i$ and $w_i$. For $j=1,2$, write $E_j$ for the number of edges in $C_j$ of colour $c^*$.
Note that (\ref{vital}) implies that $E_1-E_2\geq dr$. It is easy to see that this  implies one of $C_1$ and $C_2$ contains at least $n/r+d$ edges in the same colour,\footnote{This colour may not necessarily be $c^*$.} thereby completing the proof. \qed
\end{proofofmainthm}

\section{Concluding remarks}\label{open}
As mentioned in~\cite[Section~7]{EFLS} there are many possible directions for future research. One natural extension of our work is to seek an analogue of Theorem~\ref{mainthm} in the setting of digraphs. 
\begin{question}\label{ques1}
Given any digraph $G$ on $n$ vertices with minimum in- and outdegree at least $(1/2+1/2r+o(1))n$, and any $r$-colouring of $E(G)$, can one always ensure a Hamilton cycle in $G$ of significant colour-bias?
\end{question}
Note that the natural digraph analogues of our extremal constructions for Theorem~\ref{mainthm} show that one cannot lower the minimum degree condition in Question~\ref{ques1}.

Given an $r$-coloured $n$-vertex graph $G$ and non-negative integers $d_1,\dots,d_r$, we say that $G$ contains a \emph{$(d_1,\dots,d_r)$-coloured Hamilton cycle} if there is a Hamilton cycle in $G$ with precisely $d_i$ edges of the $i$th colour (for every $i \in [r]$). Note that the proof of Theorem~\ref{mainthm} (more precisely (\ref{vital})) ensures that given a graph $G$ as in the theorem, one can obtain at least
$dr$ distinct vectors $(d_1,\dots,d_r)$ such that $G$ has a $(d_1,\dots,d_r)$-coloured Hamilton cycle. It would be interesting to investigate this problem further. That is, given an $r$-coloured $n$-vertex graph $G$ of a given minimum degree, how many distinct vectors $(d_1,\dots,d_r)$ can we guarantee so that  $G$ contains a
$(d_1,\dots,d_r)$-coloured Hamilton cycle?

In~\cite{cpc}, the question of determining the minimum degree threshold that ensures a colour-bias $k$th power of a Hamilton cycle was raised; it would be interesting to establish whether a variant of the switching method from the proof of Theorem~\ref{mainthm} can be used to resolve this problem (for all $k \geq 2$ and $r$-colourings where $r\geq 2$).

{\bf Remark:} Since a version of this paper first appeared online, Brada\v{c}~\cite{brad}  has used the regularity method to resolve this problem asymptotically for all $k \geq 2$ when $r=2$.


\section*{Acknowledgement}
The authors are grateful to the referee for their careful review.

\medskip

{\footnotesize \obeylines \parindent=0pt
\begin{tabular}{lll}
	Andrea Freschi, Joseph Hyde \& Andrew Treglown  &\ &  Joanna Lada\\
	School of Mathematics 					&\ & Merton College		\\
	University of Birmingham				&\ & University of Oxford 	\\
	Birmingham											&\ & Oxford \\
	B15 2TT													&\ & OX1 2JD	 \\
	UK															&\ & UK
\end{tabular}
}
\begin{flushleft}
{\it{E-mail addresses}:
\tt{ $\{$axf079, jfh337, a.c.treglown$\}$@bham.ac.uk}, joanna.lada@merton.ox.ac.uk}
\end{flushleft}

\end{document}